\input amstex\documentstyle{amsppt}
\pagewidth{12.5 cm}\pageheight{19 cm}\magnification\magstep 1
\topmatter
\title Rationality properties of unipotent representations\endtitle
\author G. Lusztig\endauthor
\address Department of Mathematics, M.I.T., Cambridge, MA 02139\endaddress
\thanks{Supported by the National Science Foundation}\endthanks
\endtopmatter
\document

\define\mto{\mapsto}
\define\lra{\leftrightarrow}

\define\sm{\smallmatrix}
\define\esm{\endsmallmatrix}
\define\sub{\subset}

\define\tim{\times}
\define\ti{\tilde}
\define\nl{\newline}
\redefine\i{^{-1}}

\define\sgn{\text{\rm sgn}}
\define\tr{\text{\rm tr}}

\define\sha{\sharp}

\define\al{\alpha}
\define\bet{\beta}
\define\ga{\gamma}
\define\de{\delta}
\define\ep{\epsilon}

\define\si{\sigma}
\define\th{\theta}

\define\la{\lambda}
\define\ze{\zeta}

\define\bok{\bold k}

\define\bc{\bold C}

\define\bh{\bold H}

\define\bn{\bold N}

\define\bq{\bold Q}
\define\br{\bold R}

\define\bz{\bold Z}

\define\ca{\Cal A}

\define\cf{\Cal F}
\define\cg{\Cal G}

\define\cl{\Cal L}

\define\cu{\Cal U}

\define\cx{\Cal X}

\define\BBD{BBD}
\define\DL{DL}
\define\LC{L1}
\define\LB{L2}
\define\LP{L3}    
\define\OH{Oh} 
\define\WE{W} 
\define\ticu{\ti{\cu}_\bq}
\define\ind{\text{\rm ind}}
\head Introduction\endhead
\subhead 0.1\endsubhead
Let $\bok$ be an algebraic closure of a finite field $F_q$ with $q$ elements.
Let $G$ be a connected simple algebraic group of adjoint type over $\bok$ with
a fixed $F_q$-rational structure; let $F:G@>>>G$ be the corresponding Frobenius
map. The fixed point set $G^F$ is a finite group. Let $W$ be the Weyl group of
$G$. For $w\in W$ let $R_w$ be the character of the virtual representation 
$R^1(w)$ of $G^F$ defined in \cite{\DL, 1.5}. (The definition of $R_w$ is in
terms of $l$-adic cohomology but in fact $R_w$ has integer values and is 
independent of $l$, see \cite{\DL, 3.3}.) An irreducible representation $\rho$
of $G^F$ over $\bc$ is said to be {\it unipotent} if its character 
$\chi_\rho:G^F@>>>\bc$ occurs with $\ne 0$ multiplicity in $R_w$ for some 
$w\in W$ (see \cite{\DL, 7.8}). Let $\cu$ be the set of isomorphism classes of
unipotent representations of $G^F$. Let 
$\ticu=\{\rho\in\cu|\chi_\rho(g)\in\bq\quad\forall g\in G^F\}$. Let $\cu_\bq$ 
be the set of all $\rho\in\cu$ such that $\rho$ is defined over $\bq$ (that is,
it can be realized by a $\bq[G^F]$-module). We have $\cu_\bq\sub\ticu\sub\cu$. 

Unless otherwise specified, we assume that $G$ {\it is split over} $F_q$. The 
following is one of our results.

\proclaim{Theorem 0.2}We have $\cu_\bq=\ticu$.
\endproclaim

\subhead 0.3\endsubhead
We will also show (see 1.12) that, if $G$ is of type $A,B,C$ or $D$, then 
$\cu_\bq=\cu$. (The analogous statement is false for exceptional types.) The 
rationality of certain unipotent cuspidal representations connected with 
Coxeter elements has been proved in \cite{\LC}. The method of \cite{\LC} has 
been extended in \cite{\LP} (unpublished) to determine explicitly $\cu_\bq$ in
the general case (including non-split groups). The case where $G$ is non-split
of type $A$ has been also considered by Ohmori \cite{\OH} by another extension
of the method of \cite{\LC}. 

Our study of rationality of unipotent representations is based on the statement
that a given unipotent representation appears with multiplicity $1$ in some 
(possibly virtual) representation $R$ defined using $l$-adic cohomology and 
then using the Hasse principle. In the first method (that of \cite{\LP}), $R$ 
is a particular intersection cohomology space of a variety; see Sec.1. In the 
second method (which applies only in the cuspidal case), $R$ will be one of the
$R_w$ above; see Sec.2. In one case ($G=SO_5$ with $q$ odd), we give an 
elementary approach to rationality (without using the Hasse principle); see 
Sec.3.

\head 1. First method\endhead
\subhead 1.1\endsubhead
Let $p$ be the characteristic of $F_q$. For any prime number $l\ne p$, we 
choose an imbedding of the field $\bq_l$, the $l$-adic numbers, into $\bc$.
This allows us to regard any representation of $G^F$ over $\bq_l$ as one over 
$\bc$. Let $X$ be the flag manifold of $G$; let $F:X@>>>X$ be the map induced 
by $F:G@>>>G$. For $w\in W$ let $O_w$ be the set of all $(B,B')\in X\tim X$ 
that are in relative position $w$. As in \cite{\DL}, for any $w\in W$, let 
$X_w$ be the 
subvariety of $X$ consisting of all $B\in X$ such that $(B,F(B))\in O_w$; let
$\bar X_w$ be the closure of $X_w$ in $X$. Then $X_w,\bar X_w$ are stable under
the conjugation action of $G^F$ on $X$. Hence for any $j\in\bz$ there is an
induced action of $G^F$ on the $l$-adic cohomology with compact support 
$H^j_c(X_w,\bq_l)$ and on the $l$-adic intersection cohomology 
$\bh^j(\bar X_w,\bq_l)$. (Note that $\bar X_w$ has pure dimension $l(w)$ where
$l:W@>>>\bn$ is the length function.) Recall that $R_w$ is the character of the
virtual representation $\sum_{j\in\bz}(-1)^jH^j_c(X_w,\bq_l)$ of $G^F$.  

\proclaim{Lemma 1.2} Let $\rho\in\cu$. There exists $x\in W$ and $j\in[0,l(x)]$
such that $\rho$ appears with multiplicity $1$ in the $G^F$-module
$\bh^j(\bar X_x,\bq_l)$. 
\endproclaim
The proof is based on results of \cite{\LB}. For any $x\in W$ let $\ca_x$ be 
the virtual representation of $W$ defined in \cite{\LB, p.154,156}. For any 
virtual representation $E$ of $W$ we set $R_E=|W|\i\sum_{w\in W}\tr(w,E)R_w$.
(A $\bq$-valued class function on $G^F$.) Thus, $R_{\ca_x}$ is defined. Let 
$a:W@>>>\bn$ be as in \cite{\LB, p.178}. Assume that 

(a) $x\in W$ is such that $\rho$ appears with multiplicity $1$ in 
$(-1)^{l(x)-a(x)}R_{\ca_x}$. 
\nl
Then from \cite{\LB, 6.15, 6.17(i), 5.13(i)} we deduce that $\rho$ appears with
multiplicity $1$ in $\bh^{l(x)-a(x)}(\bar X_x,\bq_l)$. (Actually, in the 
references given, $q$ is assumed to be sufficiently large; but this assumption
is removed later in \cite{\LB}.) Thus, to prove the lemma it is enough to show
that (a) holds for some $x\in W$. Now in \cite{\LB}, the multiplicities of any
unipotent representation in $(-1)^{l(x)-a(x)}R_{\ca_x}$ have been explicitly 
described for many $x$. (See for example the tables in \cite{\LB, p.304-306} 
for types $E_8,F_4$ and the results in \cite{\LB, Ch.9} for classical types.) 
In particular, we see that (a) holds for some $x\in W$.

\proclaim{Lemma 1.3} Let $\rho\in\ticu$. Let $l$ be a prime number invertible 
in $\bok$. Let $x,j$ be as in 1.2.  

(a) $\rho$ is defined over $\bq_l$.

(b) If $j$ is even then $\rho$ is defined over $\br$. If $j$ is odd then $\rho$
is not defined over $\br$.

(c) If $j$ is even then $\rho\in\cu_\bq$.
\endproclaim
Clearly, (a) follows from 1.2. We prove (b). Let $c\in H^2(\bar X_w,\bq_l)$ be
the Chern class of an ample line bundle on $\bar X_w$ (we ignore Tate twists);
we may assume that this line bundle is the restriction of a line bundle on $X$.
Since $G^F$ acts trivially on $H^2(X,\bq_l)$ it follows that $c$ is 
$G^F$-stable. Hence the map 
$\bh^j(\bar X_x,\bq_l)@>>>\bh^{2l(x)-j}(\bar X_x,\bq_l)$ given by
$\xi\mto c^{l(x)-j}\xi$ is compatible with the $G^F$-action. This map is an
isomorphism, by the Hard Lefschetz Theorem \cite{\BBD, 5.4.10}. Let 
$(,):\bh^j(\bar X_x,\bq_l)\tim\bh^{2l(x)-j}(\bar X_x,\bq_l)$ be the Poincar\'e
duality pairing. (We again ignore Tate twists.) Then 
$\xi,\xi'\mto(\xi,c^{l(x)-j}\xi')$ is a $(-1)^j$-symmetric, non-singular,
$G^F$-invariant bilinear form
$\bh^j(\bar X_x,\bq_l)\tim\bh^j(\bar X_x,\bq_l)\to\bq_l$. This restricts to a 
$(-1)^j$-symmetric, $G^F$-invariant bilinear form on the $\rho$-isotypic part 
of $\bh^j(\bar X_x,\bq_l)$, which is non-singular, since $\rho$ is isomorphic 
to its dual (recall that $\rho\in\ticu$). This $\rho$-isotypic part is 
isomorphic to $\rho$ and (b) follows. Under the assumption of (c), we see from
(a),(b), using the Hasse principle for division algebras with centre $\bq$ that
$\rho$ is defined over $\bq$. (The Hasse principle is applicable even when
information is missing at one place, in our case at $p$-adic numbers.) The 
lemma is proved.

\proclaim{Lemma 1.4} Let $\rho\in\ticu$. Let $x,j$ be as in 1.2. Then $j$ is 
even.
\endproclaim
It is known \cite{\LB} that the parity of an integer $j$ such that $\rho$ 
appears with non-zero multiplicity in $\bh^j(\bar X_x,\bq_l)$ for some $x\in W$
is an invariant of $\rho$. Moreover, $j$ is even except if $G$ is of type $E_7$
and $\rho$ is one the two unipotent cuspidal representations or $G$ is of type
$E_8$ and $\rho$ is a component of the representation induced by one of the two
unipotent cuspidal representations of a parabolic of type $E_7$. (See 
\cite{\LB, Ch.11}.) In these exceptional cases, we have $\rho\notin\ticu$, as 
one sees using \cite{\LB, 11.2}. The lemma is proved.

\subhead 1.5\endsubhead
Now Theorem 0.2 follows immediately from 1.3(c) and 1.4. 

\subhead 1.6\endsubhead
Let $\cx$ be the set of all triples $(\cf,y,\si)$ where $\cf$ is a "family"
\cite{\LB, 4.2} of irreducible representations of $W$ (with an associated 
finite group $\cg_\cf$, see \cite{\LB, Ch.4}), $y$ is an element of $\cg_\cf$
defined up to conjugacy and $\si$ is an irreducible representation of the
centralizer of $y$ in $\cg_\cf$ defined up to isomorphism. For 
$(\cf,y,\si)\in\cx$, let $\la_{y,\si}$ be the scalar by which $y$ acts on $\si$
(a root of $1$). Let $\cx_1$ be the set of all $(\cf,y,\si)\in\cx$ such that 
$|\cf|\ne 2$ and $\la_{y,\si}=\pm 1$. If $q$ is a square, let $\cx_2$ be the 
set of all $(\cf,y,\si)\in\cx$ such that $|\cf|=2,y=1$. If $q$ is not a square,
let $\cx_2=\emptyset$. In any case, $\cx_2$ is empty unless $G$ is of type 
$E_7$ or $E_8$. Let $\cx_\bq=\cx_1\cup\cx_2$. 

In \cite{\LB, 4.23}, $\cx$ is put in a bijection 
$$(\cf,y,\si)\lra\rho_{\cf,y,\si}\tag a$$ 
with $\cu$.

\proclaim{Lemma 1.7} Assume that $\rho=\rho_{\cf,y,\si}$, 
$\rho'=\rho_{\cf',y',\si'}$ where $(\cf,y,\si)\in\cx_\bq$,
$(\cf',y',\si')\in\cx$ are distinct. Then there exists $x\in W$ such that 
$\rho,\rho'$ have different multiplicities in the $G^F$-module 
$(-1)^{l(x)-a(x)}R_{\ca_x}$. 
\endproclaim
As mentioned in the proof of 1.2, the multiplicities of various unipotent 
representations have been explicitly computed in \cite{\LB} for many $x\in W$.
From this the lemma follows easily.

\proclaim{Lemma 1.8}Let $\rho=\rho_{\cf,y,\si}$, where $(\cf,y,\si)\in\cx_\bq$.
Then $\rho\in\ticu$.
\endproclaim
Let $\ga\in\text{\rm Gal}(\bc/\bq)$. Then $\ga(\chi_\rho)=\chi_{\rho'}$ for 
some $\rho'\in\cu$. Since the character of $(-1)^{l(x)-a(x)}R_{\ca_x}$ is 
integer valued, it is fixed by $\ga$. (Here $x$ is any element of $W$.) Hence 
$\rho,\rho'$ have the same multiplicity in $(-1)^{l(x)-a(x)}R_{\ca_x}$. From 
1.7 it follows that $\rho=\rho'$. Thus, $\ga(\chi_\rho)=\chi_\rho$ for any 
$\ga\in\text{\rm Gal}(\bc/\bq)$, so that $\chi_\rho$ has rational values. The 
lemma is proved.

\proclaim{Lemma 1.9} Let $\rho=\rho_{\cf,y,\si}$, where 
$(\cf,y,\si)\notin\cx_\bq$. Then $\rho\notin\ticu$.
\endproclaim
Assume first that $\la_{y,\si}\ne\pm 1$. Then $\la_{y,\si}\notin\bq$ hence 
there exists $\ga\in\text{\rm Gal}(\bc/\bq)$ such that
$\ga(\la_{y,\si})\ne\la_{y,\si}$. Using the interpretation of $\la_{y,\si}$ 
given in \cite{\LB, 11.2}, it follows that $\ga(\chi_\rho)\ne\chi_\rho$. Hence 
$\rho\notin\ticu$. Next assume that $\la_{y,\si}=\pm 1$. Then $|\cf|=2$. 
Moreover, if $q$ is a square, then $y\ne 1$. Let $\si'$ be the character of 
$\cg_\cf=\bz/2\bz$ other than $\si$. Let $\rho'=\rho_{\cf,y,\si'}$. If 
$y\ne 1$, then by the results of \cite{\LC}, $\chi_\rho$ is carried to
$\chi_{\rho'}$ by an element of $\text{\rm Gal}(\bc/\bq)$ that takes 
$\sqrt{-q}$ to $-\sqrt{-q}$. If $y=1$, then by the known construction of 
representations of Hecke algebras in terms of $W$-graphs, $\chi_\rho$ is 
carried to $\chi_{\rho'}$ by an element of $\text{\rm Gal}(\bc/\bq)$ that takes
$\sqrt q$ to $-\sqrt q$. Hence again $\rho\notin\ticu$.

\proclaim{Proposition 1.10} Under the bijection $\cx\lra\cu$ in 1.6(a), the 
subset $\ticu$ of $\cu$ corresponds to the subset $\cx_\bq$ of $\cx$.
\endproclaim
This follows immediately from 1.9, 1.10.

Combining this proposition with 0.2, we obtain:

\proclaim{Corollary 1.11} Under the bijection $\cx\lra\cu$ in 1.6(a), the 
subset $\cu_\bq$ of $\cu$ corresponds to the subset $\cx_\bq$ of $\cx$.
\endproclaim
If $G$ is of type $A,B,C$ or $D$, then for any family $\cf$ we have 
$|\cf|\ne 2$ and the group $\cg_\cf$ is an elementary abelian $2$-group hence 
$\la_{y,\si}=\pm 1$ for any $(\cf,y,\si)\in\cx$. Thus, we have $\cx_\bq=\cx$ 
and we obtain:

\proclaim{Corollary 1.12} If $G$ is of type $A,B,C$ or $D$, then $\cu_\bq=\cu$.
\endproclaim

\subhead 1.13\endsubhead 
In this subsection we assume that $G$ is non-split. The analogues of Lemmas
1.2, 1.3 continue to hold but that of Lemma 1.4 does not. (It does in type $D$
but not in type $A$.) Also, if $G$ is non-split of type $D$, then 
$\cu_\bq=\cu$. If $G$ is non-split of type $A$ we have $\ticu=\cu$ but 
$\cu_\bq\ne\cu$ in general.

\head 2. Second method\endhead
\subhead 2.1\endsubhead
Let $n\in\bn$. Let $\la_1,\la_2,\dots,\la_a$ be a sequence of integers such
that 
$$\sum_i\la_i=n+\binom{a}{2}.\tag a$$
We define a virtual representation $[\la_1,\la_2,\dots,\la_a]$ of the symmetric
group $S_n$ as follows. If $0\le\la_1<\la_2<\dots<\la_a$, then 
$[\la_1,\la_2,\dots,\la_a]$ is the irreducible representation of $S_n$ 
corresponding to the partition $\la_1\le\la_2-1\le\dots\le\la_a-a+1$ of $n$, as
in \cite{\LB, p.81}. If $\la_1,\la_2,\dots,\la_a$ are in $\bn$ and are 
distinct, then 
$$[\la_1,\la_2,\dots,\la_a]
=\sgn(\si)[\la_{\si(1)},\la_{\si(2)},\dots,\la_{\si(a)}]$$
where $\si$ is the unique permutation of $1,2,\dots,a$ such that
$\la_{\si(1)}<\la_{\si(2)}<\dots<\la_{\si(a)}$. If $\la_1,\la_2,\dots,\la_a$ 
are not distinct, or if at least one of them is $<0$, we set 
$[\la_1,\la_2,\dots,\la_a]=0$. From the definition we see easily that
$$[\la_1,\la_2,\dots,\la_a]=[0,\la_1+1,\la_2+1,\dots,\la_a+1]$$
for any sequence of integers $\la_1,\la_2,\dots,\la_a$ such that (a) holds.

\proclaim{Lemma 2.2} Let $\la_1,\la_2,\dots,\la_a$ be a sequence of integers 
such that 2.1(a) holds. Let $w=(k)w'\in S_k\tim S_{n-k}\sub S_n$ where $(k)$ 
denotes a $k$-cycle in $S_k$ and $w'\in S_{n-k}$. We have 
$$\tr(w,[\la_1,\la_2,\dots,\la_a])=\sum_{i=1}^a
\tr(w',[\la_1,\la_2,\dots,\la_{i-1},\la_i-k,\la_{i+1},\dots,\la_a]).\tag a$$
\endproclaim
If the $\la_i$ are not distinct or if at least one of them is $<0$ then both
sides of (a) are $0$. We may assume that $0\le\la_1<\la_2<\dots<\la_a$. In this
case, (a) can be seen to be equivalent to Murnaghan's rule, see \cite{\WE}.

\subhead 2.3\endsubhead
For $n\ge 0$ let $W_n$ be the group of all permutations of 
$1,2,\dots,n,n',\dots,2',1'$ which commute with the involution $i\lra i'$ for
$i=1,\dots,n$. (We have $W_0=\{1\}$.) Given two sequences of integers 
$\la_1,\dots,\la_a$ and $\mu_1,\mu_2,\dots,\mu_b$ such that 
$$\sum_i\la_i+\sum_i\mu_i=n+\binom{a}{2}+\binom{b}{2},\tag a$$
we define a virtual representation
$$\left[\sm\la_1&\la_2&\dots&\la_a\\ \mu_1&\mu_2&\dots&\mu_b\esm\right]\tag b$$
of $W_n$ as follows. If the $\la_i$ are not distinct or if the $\mu_i$ are not
distinct or if at least one of the $\la_i$ or $\mu_i$ is $<0$ we define (b) to
be $0$. Assume now that the $\la_i\in\bn$ are distinct, and that the 
$\mu_i\in\bn$ are distinct. Then $r,\ti r$ defined by 
$$\sum_i\la_i=r+\binom{a}{2},\quad\sum_i\mu_i=\ti r+\binom{b}{2}$$
satisfy $r,\ti r\in\bn, r+\ti r=n$.
We identify $W_r\tim W_{\ti r}$ with a subgroup of $W_n$ as in 
\cite{\LB, p.82}. The virtual representation 
$[\la_1,\la_2,\dots,\la_a]\boxtimes[\mu_1,\mu_2,\dots,\mu_b]$
of $S_r\tim S_{\ti r}$ may be regarded as a virtual representation of
$W_r\tim W_{\ti r}$ via the obvious projection
$W_r\tim W_{\ti r}@>>>S_r\tim S_{\ti r}$ (see \cite{\LB, p.82}). We tensor this
with the one dimensional character of $W_r\tim W_{\ti r}$ which is the identity
on the $W_r$-factor and is the restriction of $\chi:W_n@>>>\{\pm 1\}$ (see 
\cite{\LB, p.82}) on the $W_{\ti r}$-factor. Inducing the resulting virtual 
representation from $W_r\tim W_{\ti r}$ to $W_n$, we obtain the virtual 
representation (b)
of $W_n$. Note that if $\la_1<\la_2<\dots<\la_a$ and $\mu_1<\mu_2<\dots<\mu_b$
then this is an irreducible representation; if $\si$ is a permutation of 
$1,2,\dots,a$ and $\si'$ is a permutation of $1,2,\dots,b$ then
$$\left[\sm\la_{\si(1)}&\la_{\si(2)}&\dots&\la_{\si(a)}\\ \mu_{\si'(1)}&
\mu_{\si'(2)}&\dots&\mu_{\si'(b)}\esm\right]=\sgn(\si)\sgn(\si')
\left[\sm\la_1&\la_2&\dots&\la_a\\ \mu_1&\mu_2&\dots&\mu_b\esm\right].$$ 
From the definition we see easily that
$$\left[\sm\la_1&\la_2&\dots&\la_a\\ \mu_1&\mu_2&\dots&\mu_b\esm\right]=
\left[\sm 0&\la_1+1&\la_2+1&\dots&\la_a+1\\0&\mu_1+1&\mu_2+1&\dots&\mu_b+1
\esm\right].$$
 
\proclaim{Lemma 2.4} Let $\la_1,\la_2,\dots,\la_a$ and 
$\mu_1,\mu_2,\dots,\mu_b$ be two sequences of integers such that 2.3(a) holds.
Let $w=(2k)\tim w'\in W_k\tim W_{n-k}\sub W_n$ where $0<k\le n$, $(2k)$ denotes
an element of $W_k$ whose image under the obvious imbedding $W_k\sub S_{2k}$ is
a $2k$-cycle and $w'\in W_{n-k}$ has no cycles of length $2k$ as an element of 
$S_{2n-2k}$. We have
$$\align&
\tr(w,\left[\sm\la_1&\la_2&\dots&\la_a\\ \mu_1&\mu_2&\dots&\mu_b\esm\right])
=\sum_{i=1}^a
\tr(w',\left[\sm\la_1&\la_2&\dots&\la_{i-1}&\la_i-k&\la_{i+1}&\dots\la_a
\\ \mu_1&\mu_2&\dots&\mu_b\esm\right])\\&-\sum_{i=1}^a
\tr(w',\left[\sm\la_1&\la_2&\dots&\la_a
\\ \mu_1&\mu_2&\dots&\mu_{i-1}&\mu_i-k&\mu_{i+1}&\dots&\mu_b\esm\right]).
\endalign$$
\endproclaim
This follows from Lemma 2.2, using the definitions.

\subhead 2.5\endsubhead
Let $m\in\bn$ and let $n=m^2+m$. Let $w_m\in W_n$ be an element whose image 
under the imbedding $W_n\sub S_{2n}=S_{2(m^2+m)}$ is a product of cycles 
$$(4)(8)(12)\dots(4m).$$
Let 
$$\la_1<\la_2<\dots<\la_{m+1} \text{ and } \mu_1<\mu_2<\dots<\mu_m \tag a$$
be two sequences of integers such that 
$\la_1,\la_2,\dots,\la_{m+1},\mu_1,\mu_2,\dots,\mu_m$ is a permutation of 
$0,1,2,3,\dots,2m$. Then 2.3(a) holds (with $a=m+1,b=m$ and $n=m^2+m$). 
Consider the property
$$\la_i+\la_j\ne 2m \text{ for any } i\ne j \text{ and } \mu_i+\mu_j\ne 2m
\text{ for any } i\ne j.\tag *$$

\proclaim{Lemma 2.6} In the setup of 2.5, if $(*)$ holds, then
$$\tr(w_m,\left[\sm\la_1&\la_2&\dots&\la_{m+1}\\ \mu_1&\mu_2&\dots&\mu_m\esm
\right])=(-1)^{(m^2+m)/2}.$$
If $(*)$ does not hold, then
$\tr(w_m,\left[\sm\la_1&\la_2&\dots&\la_{m+1}\\ \mu_1&\mu_2&\dots&\mu_m\esm
\right])=0.$
\endproclaim
We argue by induction on $m$. The result is clear when $m=0$. Assume now that
$m>0$. We can assume that $w_m=(4m)w_{m-1}\in W_{2m}\tim W_{n-2m}\sub W_n$ 
where $w_{m-1}\in W_{n-2m}$ is defined in a way similar to $w_m$. We apply 2.4
with $w=w_m,k=2m,w'=w_{m-1}$. Note that in the formula in 2.4, at most one term
is non-zero, namely the one in which $k=2m$ is substracted from the largest 
entry $\la_i$ or $\mu_i$ (the other terms are zero since they contain some $<0$
entry). We are in one of the four cases below.

{\it Case 1.} $2m=\la_{m+1},0=\mu_1$.

Using 2.4, we have
$$\align&
A=\tr(w_m,\left[\sm\la_1&\la_2&\dots&\la_{m+1}\\ \mu_1&\mu_2&\dots&\mu_m
\esm\right])
=\tr(w_{m-1},\left[\sm\la_1&\la_2&\dots&\la_m&0\\0&\mu_2&\dots&\mu_m\esm
\right])\\&
=(-1)^m\tr(w_{m-1},\left[\sm 0&\la_1&\la_2&\dots&\la_m\\0&\mu_2&\dots&\mu_m
\esm\right])
=(-1)^m\tr(w_{m-1},\left[\sm\la_1-1&\la_2-1&\dots&\la_m-1\\ 
\mu_2-1&\mu_3-1&\dots&\mu_m-1\esm\right]).\endalign$$
Now the induction hypothesis is applicable to 
$$\la_1-1<\la_2-1<\dots<\la_m-1 \text{ and }
\mu_2-1<\mu_3-1<\dots<\mu_m-1\tag a$$ 
instead of 2.5(a). (Clearly, 2.5(a) satisfies $(*)$ if and only if (a) 
satisfies the analogous condition). Hence, if 2.5(a) satisfies $(*)$, then 
$$A=(-1)^m(-1)^{(m^2-m)/2}=(-1)^{(m^2+m)/2}$$ 
as required. If 2.5(a) does not satisfy $(*)$, then $A=(-1)^m0=0$, as required.

{\it Case 2.} $2m=\la_{m+1},0=\la_1$.

Using 2.4, we have
$$\tr(w_m,\left[\sm\la_1&\la_2&\dots&\la_{m+1}\\ \mu_1&\mu_2&\dots&\mu_m\esm
\right])
=\tr(w_{m-1},\left[\sm 0&\la_2&\dots&\la_m&0\\ \mu_1&\mu_2&\dots&\mu_m\esm
\right])$$
and this is $0$ since $0$ appears twice in the top row. 

{\it Case 3.} $2m=\mu_m,0=\la_1$.

Using 2.4, we have
$$\align&
A=\tr(w_m,\left[\sm\la_1&\la_2&\dots&\la_{m+1}\\ \mu_1&\mu_2&\dots&\mu_m\esm
\right])=-\tr(w_{m-1},\left[\sm 0&\la_2&\dots&\la_m&\la_{m+1}\\ 
\mu_1&\mu_2&\dots&\mu_{m-1}&0\esm\right])\\&
=(-1)^m\tr(w_{m-1},\left[\sm 0&\la_2&\dots&\la_m&\la_{m+1}\\ 
0&\mu_1&\mu_2&\dots&\mu_{m-1}\esm\right])\\&
=(-1)^m\tr(w_{m-1},\left[\sm\la_2-1&\dots&\la_m-1&\la_{m+1}-1\\ 
\mu_1-1&\mu_2-1&\dots&\mu_{m-1}-1\esm\right]).\endalign$$
Now the induction hypothesis is applicable to 
$$\la_2-1<\dots<\la_m-1<\la_{m+1}-1 \text{ and }
\mu_1-1<\mu_2-1<\dots<\mu_{m-1}-1\tag b$$
instead of 2.5(a). (Clearly, 2.5(a) satisfies $(*)$ if and only if (b) 
satisfies the analogous condition.) Hence, if 2.5(a) satisfies $(*)$, then 
$$A=(-1)^m(-1)^{(m^2-m)/2}=(-1)^{(m^2+m)/2}$$
 as required. If 2.5(a) does not 
satisfy $(*)$, then $A=(-1)^m0=0$, as required.

{\it Case 4.} $2m=\mu_m,0=\mu_1$.

Using 2.4, we have
$$\tr(w_m,\left[\sm\la_1&\la_2&\dots&\la_{m+1}\\ \mu_1&\mu_2&\dots&\mu_m\esm
\right])
=-\tr(w_{m-1},\left[\sm\la_1&\la_2&\dots&\la_m&\la_{m+1}\\ 0&\mu_2&\dots&0\esm
\right])$$
and this is $0$ since $0$ appears twice in the bottom row. The lemma is proved.

\proclaim{Lemma 2.7} In the setup of 2.6, if $(*)$ holds, then 

$\sha(k\in\{1,2,\dots,m\}; \mu_k=\text{ even })=(m^2+m)/2 \mod 2$.
\endproclaim
Since (*) holds, the left hand side is equal to the number of pairs
$$(0,2m),(1,2m-1),(2,2m-2),\dots,(m-1,m+1)$$ 
in which both components are even. This equals $m/2$ if $m$ is even and 
$(m+1)/2$ if $m$ is odd. Hence it is has the same parity as $m(m+1)/2$. The 
lemma is proved.

\subhead 2.8\endsubhead
Let $m\in\bn, m\ge 1$ and let $n=m^2$. Let $w'_m\in W_n$ be an element whose 
image under the imbedding $W_n\sub S_{2n}=S_{2m^2}$ is a product of cycles 
$$(2)(6)(10)\dots(4m-2).$$
Let 
$$\la_1<\la_2<\dots<\la_m \text{ and } \mu_1<\mu_2<\dots<\mu_m\tag a$$
be two sequences of integers such that 
$\la_1,\la_2,\dots,\la_m,\mu_1,\mu_2,\dots,\mu_m$ is a permutation of 
$0,1,2,\dots,2m-1$. Then 2.3(a) holds (with $a=b=m$ and $n=m^2$). Let 
$$N=\sha(k\in\{1,2,\dots,m\}; \mu_k\ge m).\tag b$$
Consider the property
$$\la_i+\la_j\ne 2m-1 \text{ for any } i\ne j \text{ and } \mu_i+\mu_j\ne 2m-1
\text{ for any } i\ne j.\tag **$$

\proclaim{Lemma 2.9} In the setup of 2.8, if $(**)$ holds, then
$$\tr(w'_m,\left[\sm\la_1&\la_2&\dots&\la_m\\ \mu_1&\mu_2&\dots&\mu_m\esm
\right])=(-1)^{N+m(m-1)/2}.$$
If $(**)$ does not hold, then
$\tr(w'_m,\left[\sm\la_1&\la_2&\dots&\la_m\\ \mu_1&\mu_2&\dots&\mu_m\esm
\right])=0.$
\endproclaim
We argue by induction on $m$. The result is clear when $m=1$. Assume now that
$m>1$. We can assume that 
$w'_m=(4m-2)w'_{m-1}\in W_{2m-1}\tim W_{n-2m+1}\sub W_n$ 
where $w'_{m-1}\in W_{n-2m+1}$ is defined in a way similar to $w'_m$. We apply
2.4 with $w=w'_m,k=2m-1,w'=w'_{m-1}$. Note that in the formula in 2.4, at most
one term is non-zero, namely the one in which $k=2m-1$ is substracted from the
largest entry $\la_i$ or $\mu_i$ (the other terms are zero since they contain 
some $<0$ entry). We are in one of the four cases below.

{\it Case 1.} $2m-1=\la_m,0=\mu_1$.

Using 2.4, we have
$$\align&
A=\tr(w'_m,\left[\sm\la_1&\la_2&\dots&\la_m\\ \mu_1&\mu_2&\dots&\mu_m\esm
\right])
=\tr(w'_{m-1},\left[\sm\la_1&\la_2&\dots&\la_{m-1}&0\\0&\mu_2&\dots&\mu_m\esm
\right])\\&
=(-1)^{m-1}\tr(w'_{m-1},\left[\sm 0&\la_1&\la_2&\dots&\la_{m-1}\\
0&\mu_2&\dots&\mu_m\esm\right])\\&
=(-1)^{m-1}\tr(w'_{m-1},\left[\sm\la_1-1&\la_2-1&\dots&\la_{m-1}-1\\ 
\mu_2-1&\mu_3-1&\dots&\mu_m-1\esm\right]).\endalign$$
Now the induction hypothesis is applicable to 
$$\la_1-1<\la_2-1<\dots<\la_{m-1}-1 \text{ and }\mu_2-1<\mu_3-1<\dots<\mu_m-1
\tag a$$
instead of 2.8(a). (Clearly, 2.8(a) satisfies $(**)$ if and only if (a) 
satisfies the analogous condition). Let $N'$ be defined as $N$ in 2.8(b), in
terms of (a). Then $N'=N$. If 2.8(a) satisfies $(**)$, then 
$$A=(-1)^{m-1}(-1)^{(m-1)(m-2)/2}(-1)^{N'}=(-1)^{m(m-1)/2}(-1)^N$$ 
as required. If 2.8(a) does not satisfy $(**)$, then $A=(-1)^{m-1}0=0$, as 
required.

{\it Case 2.} $2m-1=\la_m,0=\la_1$.

Using 2.4, we have
$$\tr(w'_m,\left[\sm\la_1&\la_2&\dots&\la_m\\ \mu_1&\mu_2&\dots&\mu_m\esm
\right])
=\tr(w'_{m-1},\left[\sm 0&\la_2&\dots&\la_{m-1}&0\\ \mu_1&\mu_2&\dots&\mu_m
\esm\right])$$
and this is $0$ since $0$ appears twice in the top row. 

{\it Case 3.} $2m-1=\mu_m,0=\la_1$.

Using 2.4, we have
$$\align&A=
\tr(w'_m,\left[\sm\la_1&\la_2&\dots&\la_m\\ \mu_1&\mu_2&\dots&\mu_m\esm
\right])
=-\tr(w'_{m-1},\left[\sm 0&\la_2&\dots&\la_m\\ \mu_1&\mu_2&\dots&\mu_{m-1}&0
\esm\right])\\&
=(-1)^m\tr(w'_{m-1},\left[\sm 0&\la_2&\dots&\la_m\\0&\mu_1&\mu_2&\dots&
\mu_{m-1}\esm\right])\\&
=(-1)^m\tr(w'_{m-1},\left[\sm\la_2-1&\dots&\la_m-1\\ 
\mu_1-1&\mu_2-1&\dots&\mu_{m-1}-1\esm\right]).\endalign$$
Now the induction hypothesis is applicable to 
$$\la_2-1<\la_3-1<\dots<\la_m-1 \text{ and }
\mu_1-1<\mu_2-1<\dots<\mu_{m-1}-1\tag b$$
instead of 2.8(a). (Clearly, 2.8(a) satisfies $(**)$ if and only if (b) 
satisfies the analogous condition.) Let $N'$ be defined as $N$ in 2.8(b), in 
terms of (b). Then $N'=N-1$. If 2.8(a) satisfies $(**)$, then 
$$A=(-1)^m(-1)^{(m-2)(m-1)/2}(-1)^{N'}=(-1)^{m(m-1)/2}(-1)^N$$
as required. If 2.8(a) does not satisfy $(**)$, then $A=(-1)^m0=0$, as 
required.

{\it Case 4.} $2m-1=\mu_m,0=\mu_1$.

Using 2.4, we have
$$\tr(w'_m,\left[\sm\la_1&\la_2&\dots&\la_m\\ \mu_1&\mu_2&\dots&\mu_m\esm
\right])
=-\tr(w'_{m-1},\left[\sm\la_1&\la_2&\dots&\la_m\\ 0&\mu_2&\dots&0\esm
\right])$$
and this is $0$ since $0$ appears twice in the bottom row. The lemma is proved.

\proclaim{Lemma 2.10} Assume that we are in the setup of 2.9, that $(**)$ holds
and that $m=2m'$ for some integer $m'>0$. Then

(a) $\sha(k\in\{1,2,\dots,m\}; \mu_k\ge m)-
\sha(k\in\{1,2,\dots,m\}; \mu_k \text{ even })=m' \mod 2$,

(b) $\sha(k\in\{1,2,\dots,m\}; \mu_k \text{ even })=N+m(m-1)/2 \mod 2$.
\endproclaim
Among the $m'$ pairs $(0,4m'-1),(2,4m'-3),\dots,(2m'-2,2m'+1)$ there are, say,
$\al$ pairs with the first component of form $\la_i$ and second component of
form $\mu_j$ and $\bet$ pairs with the first component of form $\mu_j$ and 
second component of form $\la_i$. Clearly, $\al+\bet=m'$. Among the $m'$ pairs
$$(1,4m'-2),(3,4m'-4),\dots,(2m'-1,2m')$$ 
there are, say, $\ga$ pairs with the 
first component of form $\la_i$ and second component of form $\mu_j$ and $\de$
pairs with the first component of form $\mu_j$ and second component of form 
$\la_i$. Clearly, $\ga+\de=m'$. From the definitions we have
$$\sha(k\in\{1,2,\dots,2m'\};\mu_k\ge 2m')=\al+\ga,$$
$$\sha(k\in\{1,2,\dots,2m'\};\mu_k \text{ even })=\bet+\ga.$$
Hence the left hand side of (a) is equal to $\al+\ga-(\bet+\ga)=\al-\bet$,
which has the same parity as $\al+\beta=m'$. This proves (a). Now (b) follows
from (a) since $m'=2m'(2m'-1)/2\mod 2$. The lemma is proved.

\proclaim{Proposition 2.11} Assume that $G$ in 0.1 is of type $B_n$ or $C_n$ 
where $n=m^2+m, m\in\bn, m\ge 1$. We identify the Weyl group $W$ of $G$ with 
$W_n$ (see 2.3) in the standard way. Let $w=w_m$, see 2.5. Let $\rho$ be the 
unique unipotent cuspidal representation of $G^F$. Then $\rho$ appears with 
multiplicity $1$ in $R_w$.
\endproclaim
For any subset $J$ of cardinal $m$ of $I=\{0,1,2,\dots,2m\}$ let 
$E_J=\left[\sm\la_1&\la_2&\dots&\la_{m+1}\\ \mu_1&\mu_2&\dots&\mu_m
\esm\right]$ (an
irreducible representation of $W$) where $\mu_1<\mu_2<\dots<\mu_m$ are the 
elements of $J$ in increasing order and $\la_1<\la_2<\dots<\la_{m+1}$ are the 
elements of $I-J$ in increasing order; let $f(J)=\sha(j\in J| j \text{ even})$.
By \cite{\LB, 4.23}, the multiplicity of $\rho$ in $R_w$ is
$$2^{-m}\sum_J(-1)^{f(J)}\tr(w_m,E_J)\tag a$$
where $J$ runs over all subsets of $I$ of cardinal $m$. Using 2.6 and 2.7 we 
see that (a) equals $2^{-m}\sha(J; J\cap(2m-J)=\emptyset)=1$. The proposition 
is proved.

\proclaim{Proposition 2.12} Assume that $G$ in 0.1 is of type $D_n$ where 
$n=m^2, m=2m', m'\in\bn, m'\ge 1$. We identify the Weyl group $W$ of $G$ with 
the subgroup of $W_n$ consisting of all permutations  $w\in W_n$ such that
$$\sha(k\in\{1,2,\dots,n\};w(k)\in\{1',2',\dots,n'\})$$
is even. (A subgroup of index $2$.) Let $w=w'_m$, see 2.8. (We have 
$w'_m\in W$.) Let $\rho$ be the unique unipotent cuspidal representation of 
$G^F$. Then $\rho$ appears with multiplicity $1$ in $R_w$.
\endproclaim
For any subset $J$ of cardinal $m$ of $I=\{0,1,2,\dots,2m-1\}$ let $E_J$ be the
restriction of $\left[\sm\la_1&\la_2&\dots&\la_{m+1}\\ \mu_1&\mu_2&\dots&\mu_m
\esm\right]$
from $W_n$ to $W$ (an irreducible representation of $W$) where 
$\mu_1<\mu_2<\dots<\mu_m$ are the elements of $J$ in increasing order and 
$\la_1<\la_2<\dots<\la_m$ are the elements of $I-J$ in increasing order; let 
$f(J)=\sha(j\in J| j \text{ even})$. Note that $E_J=E_{I-J}$ and $f(J)=f(I-J)$.
By \cite{\LB, 4.23}, the multiplicity of $\rho$ in $R_w$ is
$$2^{-m}\sum_J(-1)^{f(J)}\tr(w_m,E_J)\tag a$$
where $J$ runs over all subsets of $I$ of cardinal $m$. Using 2.9 and 2.10 we 
see that (a) equals $2^{-m}\sha(J; J\cap(2m-J)=\emptyset)=1$. The proposition 
is proved.

We return to the general case.

\proclaim{Theorem 2.13} Let $\rho$ be a unipotent cuspidal representation of
$G^F$. There exists $w\in W$ such that $\rho$ appears with multiplicity $1$ in
$R_w$.
\endproclaim
In view of 2.11, 2.12, we may assume that $G$ is of exceptional type. If 
$\rho\in\cu$ is cuspidal, we have $\rho_{\cf,y,\si}$ where $\cf$ is independent
of $\rho$; we will write $\rho_{y,\si}$ instead of $\rho_{\cf,y,\si}$; for the
pairs $(y,\si)$ we will use the notation of \cite{\LB, 4.3}. For $w\in W$ we 
denote by $|w|$ the characteristic polynomial of $w$ in the reflection 
representation of $W$. (A product of cyclotomic polynomials $\Phi_d$.) In the
cases that appear below, $|w|$ determines uniquely $w$ up to conjugacy.

\medpagebreak

{\it Type $E_6$.}

$|w|=\Phi_{12}\Phi_3$:  $\rho_{g_3,\th^{\pm 1}}$.

\medpagebreak

{\it Type $E_7$.}

$|w|=\Phi_{18}\Phi_2$:  $\rho_{g_2,1},\rho_{g_2,\ep}$.

\medpagebreak

{\it Type $E_8$.}

$|w|=\Phi_{30}$:  $\rho_{g_5,\ze^j},j=1,2,3,4$; $\rho_{g_6,\th^{\pm 1}}$.

$|w|=\Phi_{24}$:  $\rho_{g_4,i^{\pm 1}}$.

$|w|=\Phi_{18}\Phi_6$:  $\rho_{g_3,\ep\th^{\pm 1}}$.

$|w|=\Phi_{12}^2$:  $\rho_{g'_2,\ep}$.

$|w|=\Phi_{12}\Phi_6^2$:  $\rho_{g_2,-\ep}$.

$|w|=\Phi_6^4$:  $\rho_{1,\la^4}$.

\medpagebreak

{\it Type $F_4$.}

$|w|=\Phi_{12}$:  $\rho_{g_3,\th^{\pm 1}}$, $\rho_{g_4,i^{\pm 1}}$.

$|w|=\Phi_6^2$:  $\rho_{g'_2,\ep}$.

$|w|=\Phi_8$:  $\rho_{g_2,\ep}$.

$|w|=\Phi_4^2$:  $\rho_{1,\la^3}$.

\medpagebreak

{\it Type $G_2$.}

$|w|=\Phi_6$:  $\rho_{g_3,\th^{\pm 1}}$, $\rho_{g_2,\ep}$.

$|w|=\Phi_3$:  $\rho_{1,\la^2}$.

\medpagebreak

In each case, one can compute the multiplicity of $\rho_{y,\si}$ in $R_w$ for
$w$ in the same row, using \cite{\LB, 4.23}; for the computation we need the
character table of $W$ and the explicit entries of the non-abelian Fourier 
transform \cite{\LB, p.110-113}. The result in each case is $1$. This completes
the proof.

\proclaim{Theorem 2.14} (a) Assume that $\rho\in\ticu$ is cuspidal. Then
$\rho\in\cu_\bq$.

(b) If $G$ is of type $B,C$ or $D$ and $\rho$ is a unipotent cuspidal
representation then $\rho\in\cu_\bq$.
\endproclaim
We prove (a). Let $w\in W$ be such that $\rho$ appears with multiplicity $1$ in
$R_w$ (see 2.13). Since $R_w$ is the character of a virtual representation 
defined over $\bq_l$, it follows that $\rho$ is defined over $\bq_l$. (Here $l$
is any prime $\ne p$.) Using the Hasse principle it is then enough to show that
$\rho$ is defined over $\br$. The following argument is inspired by an argument
of Ohmori \cite{\OH}. Let $x$ be an element of minimal length in $W$ such that
$\rho$ appears with odd multiplicity in $R_x$. (Such $x$ exists by 2.13.) The
multiplicity of $\rho$ in $\sum_j(-1)^j\bh^j(\bar X_x,\bq_l)$ is equal to the
multiplicity of $\rho$ in $\sum_j(-1)^jH^j_c(X_x,\bq_l)$ plus an integer linear
combination of the multiplicities of $\rho$ in 
$\sum_j(-1)^jH^j_c(X_{x'},\bq_l)$ for various $x'$ of strictly smaller length 
than $x$ (these multiplicities are even, by the choice of $x$). It follows that
the multiplicity of $\rho$ in $\sum_j(-1)^j\bh^j(\bar X_x,\bq_l)$ is odd. 
Recall that $\bar X_x$ has pure dimension $l(x)$. By Poincar\'e duality, the 
$G^F$-modules $\bh^j(\bar X_x,\bq_l),\bh^{2l(x)-j}(\bar X_x,\bq_l)$ are dual to
each other. Hence they contain $\rho$ with the same multiplicity (recall that 
$\rho$ is self-dual.) It follows that the multiplicity of $\rho$ in 
$\bh^{l(x)}(\bar X_x,\bq_l)$ is odd. Now $\bh^{l(x)}(\bar X_x,\bq_l)$ admits a
$(-1)^{l(x)}$-symmetric, non-degenerate, $G^F$-invariant $\bq_l$-bilinear form
with values in $\bq_l$. (Actually, $l(x)$ is even, by 1.4.) Since $\rho$ is 
self-dual and has odd multiplicity in $\bh^{l(x)}(\bar X_x,\bq_l)$, an argument
in \cite{\OH} shows that $\rho$ itself (regarded as a $\bq_l[G^F]$-module) 
admits a symmetric, non-degenerate, $G^F$-invariant $\bq_l$-bilinear form with
values in $\bq_l$. It follows that $\rho$ is defined over $\br$. This proves 
(a).

In the setup of (b), $\rho$ is unique up to isomorphism hence it is
automatically in $\ticu$. Thus, (a) is applicable and $\rho\in\cu_\bq$. The
theorem is proved.

\subhead 2.15\endsubhead
In this subsection we assume that $G$ is non-split. The analogue of 2.13
continues to hold for $G$. But the analogue of 2.14(a) fails if $G$ is 
non-split of type $A$. 

\subhead 2.16\endsubhead
A statement like 2.13 was made without proof in \cite{\LB, p.356} (for not 
necessarily split $G$). In that statement, the assumption that $\rho$ is 
cuspidal was missing. That assumption is in fact necessary, as 2.17(ii) below
(for $G$ of type $C_4$) shows.

\proclaim{Lemma 2.17}(i) Let $\ep:W_2\tim W_2@>>>\{\pm 1\}$ be a character. 
Then 
$$\tr(w,\ind_{W_2\tim W_2}^{W_4}(\ep))\in 2\bz\tag a$$
for all $w\in W_4$.

(ii) Let $E=\left[\sm 1&2\\2&{}\esm\right]$ (an ireducible representation of 
$W_4$). Then $R_E$ is of the form $\chi_\rho$ for some $\rho\in\cu$. The 
multiplicity of $\rho$ in $R_w$ is even for any $w\in W=W_4$.
\endproclaim
The residue class $\mod 2$ of the left hand side of (a) is clearly independent
of the choice of $\ep$. Hence to prove (a) we may assume that $\ep=1$. Let 
$\pi:W_4@>>>S_4$ be the canonical homomorphism. We have
$\tr(w,\ind_{W_2\tim W_2}^{W_4}(1))=\tr(\pi(w),\ind_{S_2\tim S_2}^{S_4}(1))$.
But if $y\in S_4$, then $\tr(y,\ind_{S_2\tim S_2}^{S_4}(1))$ is $6$ if $y=1$, 
is $2$ if $y$ has order $2$ and is $0$, otherwise; in particular, it is even 
for any $y$. This proves (i).

In (ii), the multiplicity of $\rho$ in $R_w$ is $\tr(w,E)$ that is, the left 
hand side of (a) for a suitable $\ep$. Hence it is even by (i). The lemma is
proved.

\head 3. An example in $SO_5$\endhead
\subhead 3.1\endsubhead
In this section we assume that $p\ne 2$ and that $G=SO(V)$ where $V$ is a
$5$-dimensional $\bok$-vector space with a fixed $F_q$-rational structure and a
fixed non-degenerate symmetric bilinear form $(,)$ defined over $F_q$. Let $C$
be the set of all $g\in G$ such that $g=su=us$ where $-s\in O(V)$ is a 
reflection and $u\in SO(V)$ has Jordan blocks of sizes $2,2,1$. Then $C$ is a 
conjugacy class in $G$ and $F(C)=C$. A line $L$ in $V(F_q)$ is said to be of 
type $1$ if $(x,x)\in F_q^2-0$ for any $x\in L-\{0\}$ and of type $-1$ if 
$(x,x)\in F_q-F_q^2$ for any $x\in L-\{0\}$. Let $\cl_1$ (resp. $\cl_{-1}$) be 
the set of lines of type $1$ (resp. $-1$) in $V(F_q)$. For 
$\ep,\de\in\{1,-1\}$, let $C^{\ep,\de}$ be the set of all $g\in C^F$ such that
the line $L$ in $V(F_q)$ such that $g|_L=1$ is in $\cl_\ep$ and any line $L$ in
$V(F_q)$ such that $g|_L=-1$ and $(L,L)\ne 0$ is in $\cl_\de$. Then 
$C^{\ep,\de}$ is a conjugacy class of $G^F$ and $C^F$ is union of the four 
conjugacy classes $C^{1,1},C^{1,-1},C^{-1,1},C^{-1,-1}$. We define a class 
function $\phi:G^F@>>>\bz$ by $\phi(g)=2\de q$ if $g\in C^{\ep,\de}$ and 
$\phi(g)=0$ if $g\in G-C^F$. (This is the characteristic function of a cuspidal
character sheaf on $G$.)

Let $O_+$ (resp. $O_-$) be the stabilizer in $G$ of a $4$-dimensional subspace
of $V$ defined over $F_q$ on which $(,)$ is non-degenerate and split (resp.
non-split). Let $\det:O_+@>>>\{\pm 1\}$ (resp. $\det:O_-@>>>\{\pm 1\}$) be the
unique nontrivial homomorphism of algebraic groups. The restriction of $\det$ 
to $O_+^F$ or $O_-^F$ is denoted again by 
$\det$. Consider the virtual representation 
$$\Phi=\ind_{O_+^F}^{G^F}(1)-\ind_{O_+^F}^{G^F}(\det)-\ind_{O_-^F}^{G^F}(1)
+\ind_{O_-^F}^{G^F}(\det)$$
of $G^F$. For $g\in G^F$ we have
$$\tr(g,\Phi)=2\sha(L\in\cl_1; g|_L=-1\}-2\sha(L\in\cl_{-1}; g|_L=-1\}$$. 
It follows easily that $\tr(g,\Phi)=\phi(g)$.

Let $\th$ be the unique unipotent cuspidal representation of $G^F$. Then $\th$
appears with multiplicity $1$ in $\phi$. It follows that $\th$ appears with 
multiplicity $1$ in $\Phi$. Since the character of $\th$ is $\bq$-valued 
and $\Phi$ is a difference of two representations defined over $\bq$, it 
follows that $\th$ {\it is defined over} $\bq$. Thus we have proved the 
rationality of $\th$ without using the Hasse principle.

\subhead 3.2\endsubhead
Assume now that $q=3$. Then $SO(V)$ is isomorphic to a Weyl group $W$ of type 
$E_6$ while $O_+^F$ is isomorphic to a Weyl group of type $F_4$ and $O_-^F$ is 
isomorphic to a Weyl group of type $A_5\tim A_1$ (imbedded in the standard way
in the $W$). Now $\th$ corresponds to the $6$-dimensional reflection 
representation of $W$ (Kneser). Its restriction to the Weyl group of type $F_4$
contains no one dimensional invariant subspace while its restriction to the 
Weyl group of type $A_5\tim A_1$ splits into a $5$-dimensional irreducible
representation and a non-trivial $1$ dimensional representation. Since $\th$ 
has multiplicity $1$ in $\Phi$ (see 3.1) it follows that $\th$ has multiplicity
$1$ in $\ind_{O_-^F}^{G^F}(\det)$.


\Refs
\widestnumber\key{BBD}
\ref\key\BBD\by A. A. Beilinson, J. Bernstein and P. Deligne\paper Faisceaux
pervers\jour Asterisque\vol 100\yr 1982\endref
\ref\key\DL\by P. Deligne and G. Lusztig\paper Representations of reductive 
groups over finite fields\jour Ann. Math\vol 103\yr 1976\pages 103-161\endref
\ref\key\LC\by G. Lusztig\paper Coxeter orbits and eigenspaces of Frobenius
\jour Inv. Math.\vol 38\yr 1976\pages 101-159\endref
\ref\key\LB\by G. Lusztig\book Characters of reductive groups over a finite
field, Ann. Math. Studies 107\publ Princeton Univ. Press\yr 1984\endref
\ref\key\LP\by G. Lusztig\paper lecture at the U.S.-France Conference on
Representation Theory, Paris 1982\finalinfo unpublished\endref
\ref\key\OH\by Z. Ohmori\paper The Schur indices of the cuspidal unipotent
characters of the finite unitary groups\jour Proc. Japan Acad. A (Math. Sci.)
\vol 72\yr 1996\pages 111-113\endref
\ref\key\WE\by H. Weyl\book The Classical Groups\publ Princeton Univ. Press
\endref
\endRefs
\enddocument